\documentclass[12pt,letterpaper,oneside,reqno]{amsart}
\usepackage{amsfonts}
\usepackage{amsmath}
\usepackage{amssymb}
\usepackage{amsthm}
\usepackage{float}
\usepackage{mathrsfs}
\usepackage{colonequals}
\usepackage[font=small,labelfont=bf]{caption}
\usepackage[left=1in,right=1in,bottom=1in,top=1in]{geometry}
\usepackage[pdfpagelabels,hyperindex,colorlinks=true,linkcolor=blue,urlcolor=magenta,citecolor=green]{hyperref}
\newcommand \bernoulli [2][B] {{#1}\sb{#2}}

\newcommand \coeffA [3][A] {{\mathbf{#1}} \sb{#2,#3}}

\newtheorem{thm}{Theorem}[section]

\newtheorem{conj}[thm]{Conjecture}

\numberwithin{equation}{section}
\title[An unusual identity for odd-powers]
{An unusual identity for odd-powers}
\author[Petro Kolosov]{Petro Kolosov}
\email{kolosovp94@gmail.com}
\urladdr{https://kolosovpetro.github.io}
\date{\today}
\hypersetup{
pdftitle={An unusual identity for odd-powers},
pdfsubject={Discrete Mathematics, Number Theory, Combinatorics},
pdfauthor={Petro Kolosov},
pdfkeywords={
Binomial theorem,
Power function,
Polynomials,
Polynomial identities,
Multinomial theorem,
Binomial coefficient,
Bernoulli number,
Pascal's triangle,
Faulhaber's formula,
Power sums,
Worpitzky identity,
Binomial expansion
}
}
\begin{document}
    \begin{abstract}
        In this manuscript we provide a new polynomial pattern.
        This pattern allows to find a polynomial expansion of the form
        \[x^{2m+1} = \sum_{k=1}^{x}\sum_{r=0}^{m} \mathbf{A}_{m,r} k^r (x-k)^r,\]
        where $x,m\in\mathbb{N}$ and $\mathbf{A}_{m,r}$ is real coefficient.
    \end{abstract}

    \maketitle

    \section{Introduction and Main Results} \label{sec:introduction}
    We begin our mathematical journey from investigation of the pattern in terms of finite differences
    $\Delta$ of cubes $x^3$.
    Consider the table of finite differences $\Delta$ of the polynomial $x^3$
    \begin{table}[H]
        \begin{tabular}{c|cccc}
            $x$ & $x^3$ & $\Delta(x^3)$ & $\Delta^2(x^3)$ & $\Delta^3(x^3)$ \\ [3px]
            \hline
            0   & 0     & 1             & 6               & 6               \\
            1   & 1     & 7             & 12              & 6               \\
            2   & 8     & 19            & 18              & 6               \\
            3   & 27    & 37            & 24              & 6               \\
            4   & 64    & 61            & 30              & 6               \\
            5   & 125   & 91            & 36              &               \\
            6   & 216   & 127           &                 &                \\
            7   & 343   &               &                 &
        \end{tabular}
        \caption{Table of finite differences $\Delta$ of $x^3$} \label{tab:table}
    \end{table}
    It is easy to observe that finite differences $\Delta$ of polynomial $x^3$ may be expressed according
    to the pattern
    \begin{align*}
        \Delta(0^3) &= 1+6 \cdot 0 \\
        \Delta(1^3) &= 1+6\cdot0+6\cdot1 \\
        \Delta(2^3) &= 1+6\cdot0+6\cdot1+6\cdot2 \\
        \Delta(3^3) &= 1+6\cdot0+6\cdot1+6\cdot2+6\cdot3 \\
        &\; \; \vdots \\
        \Delta(x^3) &= 1+6\cdot0+6\cdot1+6\cdot2+6\cdot3+\cdots+6\cdot x
    \end{align*}
    Furthermore, the polynomial $x^3$ turns into
    \begin{align*}
        x^3&=(1+6\cdot0)+(1+6\cdot0+6\cdot1)+(1+6\cdot0+6\cdot1+6\cdot2)+\cdots \\
        &+(1+6\cdot0+6\cdot1+6\cdot2+\cdots+6\cdot(x-1))
    \end{align*}
    If we compact above expression, we get
    \begin{align*}
        x^3&=x+(x-0)\cdot6\cdot0+(x-1)\cdot6\cdot1+(x-2)\cdot6\cdot2+\cdots&+(x-(x-1))\cdot6\cdot(x-1)
    \end{align*}
    Therefore, we can consider $x^3$ as
    \begin{equation*}
        x^3 = \sum_{k=0}^{x-1} 6k(x-k) + 1 = \sum_{k=1}^{x} 6k(x-k) + 1
    \end{equation*}
    since that term $k(x-k)$ is symmetrical over $x$ for $k=0,1,..,x$.
    Now we can assume that $\sum_{k} 6k(x-k) + 1$ has the implicit form as
    \begin{equation*}
        x^3 = \sum_{k} \coeffA{1}{1} k^1(x-k)^1 + \coeffA{1}{0} k^0(x-k)^0,
    \end{equation*}
    where $\coeffA{1}{1} = 6, \coeffA{1}{0} = 1$.
    The main problem we meet is to generalize above pattern for some power $t>3$.
    Let be a conjecture
    \begin{conj}
        For every $m\in\mathbb{N}$ there are exist $\coeffA{m}{0}, \coeffA{m}{1},..., \coeffA{m}{m}$ such that
        \begin{equation*}
            \label{conjecture}
            x^{2m+1} = \sum_{k=1}^{x} \coeffA{m}{0} k^0 (x-k)^0 + \coeffA{m}{1} (x-k)^1 + \coeffA{m}{2} k^2 (x-k)^2
            + \cdots + \coeffA{m}{m} k^m (x-k)^m.
        \end{equation*}
    \end{conj}
    Consider the case $m=1$
    \[
        x^3 = \sum_{k=1}^{x} \coeffA{1}{1} k^1(x-k)^1 + \coeffA{1}{0} k^0(x-k)^0
    \]
    We evaluate the coefficients $\coeffA{1}{0}, \coeffA{1}{1}$ as follows
    \begin{align*}
        x^3 &= \sum_{k=1}^{x} \coeffA{1}{1} kx - \coeffA{1}{1} k^2 + \coeffA{1}{0}\\
        x^3 &= \coeffA{1}{1}x \sum_{k=1}^{x} k - \coeffA{1}{1} \sum_{k=1}^{x} k^2 + \sum_{k=1}^{x} \coeffA{1}{0}
    \end{align*}
    Furthermore, by means of Faulhaber's formula~\cite{kargin2020formulas} we collapse the sums
    \begin{align*}
        x^3 &= \coeffA{1}{1} x \frac{x^2+x}{2} - \coeffA{1}{1} \frac{2x^3+3x^2+x}{6} + \coeffA{1}{0}x\\
        x^3 &= \coeffA{1}{1} \frac{3x^3+3x^2}{6} - \coeffA{1}{1} \frac{2x^3 + 3x^2 +x}{6} + \coeffA{1}{0}x\\
        x^3 &= \coeffA{1}{1} \frac{x^3-x}{6} + \coeffA{1}{0}x
    \end{align*}
    Multiply both part by 6 and moving $6x^3$ to the left part gives
    \begin{align*}
        &\coeffA{1}{1} x^3 - \coeffA{1}{1} x + 6\coeffA{1}{0}x - 6x^3 = 0\\
        &x^3(\coeffA{1}{1} - 6) + x(6\coeffA{1}{0} - \coeffA{1}{1}) = 0
    \end{align*}
    Since that $x\geq 1$ we have to solve the following system of equations
    \[
        \begin{cases}
            \coeffA{1}{1} - 6 = 0\\
            6\coeffA{1}{0} - \coeffA{1}{1} = 0
        \end{cases}
    \]
    Which gives $\coeffA{1}{1} = 6$ and $\coeffA{1}{0} = 1$.
    Therefore,
    \[
        x^{3} = \sum_{k=1}^{x} 6k(x-k) + 1.
    \]
    Consider the case $m=2$.
    Let be
    \[
        x^5 = \sum_{k=1}^x \coeffA{2}{2} k^2(x-k)^2 + \coeffA{2}{1} k(x-k) + \coeffA{2}{0}
    \]
    As above, we replace the sums by means of Faulhaber's formula~\cite{kargin2020formulas}
    \begin{gather*}
        \frac{\coeffA{2}{2} x^5 - \coeffA{2}{2} x + 30\coeffA{2}{0} x }{30} +
        \frac{\coeffA{2}{1} x^3 - \coeffA{2}{1} x}{6} - x^5 = 0\\
        \coeffA{2}{2} x^5 - \coeffA{2}{2} x + 30\coeffA{2}{0} x + 5\coeffA{2}{1} x^3 - 5\coeffA{2}{1} x - 30x^5 = 0
    \end{gather*}
    Substituting $x = 1$ we get $30\coeffA{2}{0} - 30 = 0$, hence $\coeffA{2}{0} = 1$.
    Moving $x$ out of the braces we get
    \[
        x^5 (\coeffA{2}{2} - 30) + 5\coeffA{2}{1}x^3 - x(\coeffA{2}{2} - 30\coeffA{2}{0} + 5\coeffA{2}{1}) = 0
    \]
    It produces the following system of equations
    \[
        \begin{cases}
            \coeffA{2}{2} - 30 = 0\\
            \coeffA{2}{2} - 30\coeffA{2}{0} + 5\coeffA{2}{1} = 0
        \end{cases}
    \]
    Which leads to the conclusion $\coeffA{2}{2} = 30, \coeffA{2}{1} = 0, \coeffA{2}{0} = 1$.
    Finally, we get another polynomial identity
    \[
        x^5 = \sum_{k=1}^{x} 30k^2(x-k)^2 + 1.
    \]
    \begin{thm}
        For every $x,m\in\mathbb{N}$ there are $\coeffA{m}{0}, \coeffA{m}{1},\dots,\coeffA{m}{m}$,
        such that
        \begin{equation*}
            x^{2m+1} = \sum_{k=1}^{x}\sum_{r=0}^{m} \coeffA{m}{r} k^r (x-k)^r,
            \label{eq:theorem}
        \end{equation*}
        where $\coeffA{m}{r}$ is real coefficient.
    \end{thm}
    Therefore, conjecture~\ref{conjecture} is true.
    For $m > 0$ we have the following identities
    \begin{align*}
        x^3 &= \sum_{k=1}^{x} 6k(x-k) + 1 \\
        x^5 &= \sum_{k=1}^{x} 30k^2(x-k)^2 + 1 \\
        x^7 &= \sum_{k=1}^{x} 140 k^3 (x-k)^3 - 14k(x-k) + 1 \\
        x^9 &= \sum_{k=1}^{x} 630 k^4(x-k)^4 - 120k(x-k) + 1 \\
        x^{11} &= \sum_{k=1}^{x} 2772 k^5(x-k)^5 + 660 k^2(x-k)^2 - 1386k(x-k) + 1 \\
        x^{13} &= \sum_{k=1}^{x} 51480 k^5(x-k)^7 - 60060 k^3(x-k)^3 + 491400k^2(x-k)^{2} - 450054k(x-k) + 1 \\
    \end{align*}
    Moreover, since that $k(x-k)$ is symmetric over $x$, we can conclude that
    \[
        x^{2m+1} = \sum_{k=1}^{x}\sum_{r=0}^{m} \coeffA{m}{r} k^r (x-k)^r
        = \sum_{k=0}^{x-1}\sum_{r=0}^{m} \coeffA{m}{r} k^r (x-k)^r
    \]
    Coefficients $\coeffA{m}{r}$ may be calculated recursively~\cite{kolosov2016link} as follows
    \begin{equation}
        \label{eq:def_coeff_a}
        \coeffA{m}{r} \colonequals
        \begin{cases}
            (2r+1) \binom{2r}{r}, & \text{if } r=m; \\
            (2r+1) \binom{2r}{r} \sum_{d=2r+1}^{m} \coeffA{m}{d} \binom{d}{2r+1} \frac{(-1)^{d-1}}{d-r}
            \bernoulli{2d-2r}, & \text{if } 0 \leq r<m; \\
            0, & \text{if } r<0 \text{ or } r>m,
        \end{cases}
    \end{equation}
    where $\bernoulli{t}$ are Bernoulli numbers~\cite{WeissteinBernoulli}.
    It is assumed that $\bernoulli{1}=\frac{1}{2}$.
    Reader may found more information concerning coefficients $\coeffA{m}{r}$
    in OEIS~\cite{SloaneOeis302971,SloaneOeisA304042}.
    To check formulas, use the Wolfram mathematica Package\cite{GitMathematica}.
    \bibliographystyle{unsrt}
    \bibliography{an_unusual_identity_refs}

\begin{thebibliography}{1}

\bibitem{kargin2020formulas}
Levent Kargın, Ayhan Dil, and Mümün Can.
\newblock {Formulas for sums of powers of integers and their reciprocals}.
\newblock {\em arXiv preprint arXiv:2006.01132}, 2020.
\newblock \url {https://arxiv.org/abs/2006.01132}.

\bibitem{kolosov2016link}
Petro Kolosov.
\newblock {On the link between Binomial Theorem and Discrete Convolution of
  Polynomials}.
\newblock {\em arXiv preprint arXiv:1603.02468}, 2016.
\newblock \url {https://arxiv.org/abs/1603.02468}.

\bibitem{WeissteinBernoulli}
Eric~W Weisstein.
\newblock {"Bernoulli Number." From MathWorld -- A Wolfram Web Resource.}
\newblock \url {http://mathworld.wolfram.com/BernoulliNumber.html}.

\bibitem{SloaneOeis302971}
OEIS Foundation~Inc. (2020).
\newblock {The On-Line Encyclopedia of Integer Sequences}.
\newblock \url {https://oeis.org/A302971}.

\bibitem{SloaneOeisA304042}
OEIS Foundation~Inc. (2020).
\newblock {The On-Line Encyclopedia of Integer Sequences}.
\newblock \url {https://oeis.org/A304042}.

\bibitem{GitMathematica}
Petro Kolosov.
\newblock {Supplementary Mathematica Programs}.
\newblock 2020.
\newblock \url {https://github.com/kolosovpetro/Mathematica-scripts}.

\end{thebibliography}

\end{document}